\documentstyle[amsfonts,11pt]{article}
\newcommand{\R}{{\mathbb{R}}}
\newcommand{\C}{{\mathbb{C}}}
\newcommand{\F}{{\mathbb{F}}}
\newcommand{\g}{{\mathfrak g}}

\newcommand{\h}{{\mathfrak h}}

\newcommand{\fg}{{\mathfrak g}}
\newcommand{\fh}{{\mathfrak h}}

\newcommand{\fso}{{\mathfrak s}{\mathfrak o}}

\newcommand{\ad}{{\mathrm a}{\mathrm d\,}}

\newcommand{\Char}{{\mathrm C}{\mathrm h}{\mathrm a}{\mathrm r\,}}

\newcommand{\Ker}{{\mathrm K}{\mathrm e}{\mathrm r\,}}

\newcommand{\pf}{{\noindent\it Proof.\,\,\, }}
\newcommand{\qed}{\hfill{$\Box$}}
 \makeatother

\newtheorem{thm}{Theorem}[section]
\newtheorem{Def}[thm]{Definition}

\newtheorem{prop}[thm]{Proposition}

\newtheorem{lem}[thm]{Lemma}
\newtheorem{rem}[thm]{Remark}
\newtheorem{cor}[thm]{Corollary}

\begin{document}

\title{\Large{\bf {On Pseudo-Riemannian Lie Algebras:
\\ A Class of New Lie-admissible Algebras}}
\author{{Chen Zhiqi\thanks{E-mail: chenzhiqi@nankai.edu.cn.} } and
{Zhu Fuhai\thanks {E-mail: zhufuhai@nankai.edu.cn. Corresponding
author. This work is partially supported by (No. 10501025) the
National Science Foundation of China.}}\\{\footnotesize School of
Mathematical Sciences and LPMC, Nankai University, Tianjin 300071,
China}} }

\date{} \maketitle
\abstract{M. Boucetta introduced the notion of pseudo-Riemannian
Lie algebra in~\cite{MB2} when he studied the line Poisson
structure on the dual of a Lie algebra. In this paper, we redefine
pseudo-Riemannian Lie algebra, which, in essence, is a class of
new Lie-admissible algebras and prove that all pseudo-Riemannian
Lie algebras are solvable. Using our main result and method, we
prove some of M. Boucetta's results in~\cite{MB2,MB3} in a new
simpler way, and we give an explicit construction of Riemann-Lie
algebras.



\date{} \maketitle
\section{Introduction}

M. Boucetta introduced the notion of Poisson manifold with
compatible pseudo-metric in \cite{MB1} and a new class of Lie
algebras called pseudo-Riemannian Lie algebras in \cite{MB2}. He
proved that a linear Poisson structure on the dual of a Lie
algebra has a compatible pseudo-metric if and only if the Lie
algebra is a pseudo-Riemannian Lie algebra, and that the Lie
algebra obtained by linearizing at a point in a Poisson manifold
with compatible pseudo-metric is a pseudo-Riemannian Lie algebra.
See \cite{MB2} for more details of pseudo-Riemannian Poisson
manifolds and their relationship with pseudo-Riemannian Lie
algebras. Furthermore, he proved in \cite{MB3} five equivalent
conditions for $\g$ being a Riemann-Lie algebra. In this paper we
prove that any pseudo-Riemannian Lie algebra is solvable and use
our method to give a new simple proof of M. Boucetta's result.

The paper is organized as follows. In section~\ref{secp}, we
collect some basic definitions of pseudo-Riemannian Lie algebra
and translate it into our language, which is easier to deal with.
In our new definition, pseudo-Riemannian Lie algebra is a class of
Lie-admissible algebras with a compatible nondegenerate symmetric
bilinear form. In section~\ref{secm}, we prove that there is no
pseudo-Riemannian Lie algebra structure on any semi-simple Lie
algebra, which implies our main result via Levi decomposition.
In~\cite{MB3}, Boucetta classified Riemann-Lie algebras, which we
discuss in section~\ref{sec4}. Boucetta claimed the classification
of pseudo-Riemannian Lie algebras of dimension 2 and 3 without
proof in~\cite{MB2} since his proof is long, while using our
method, we give an explicit classification and concrete formula
for these pseudo-Riemannian Lie algebras. Actually,
Theorem~1.6~2b) in~\cite{MB2}, which classifies 3-dimensional
pseudo-Riemannian Lie algebras, is not quite correct.

\section{Preliminary}\label{secp}

Let $\fg$ be a Lie algebra over a field ${\mathbb F}$ and let $(\
,\ )$ be a non-degenerate symmetric bilinear form on $\fg$. If
$\Char{\mathbb F}\neq 2$, define a bilinear map
$A:\fg\times\fg\rightarrow\fg$ by
\begin{equation}\label{d1}
2(A_uv,w)=([u,v],w)+([w,u],v)+([w,v],u),\ \forall u,v,w\in\g
\end{equation}
One can easily prove that the above equality is equivalent to the
following:
\begin{equation}\label{P1}
A_uv-A_vu=[u,v],
\end{equation}
\begin{equation}\label{P2}
(A_uv,w)+(v,A_uw)=0.
\end{equation}

\begin{Def}
A pair $(\fg,(\ ,\ ))$ (or $\g$ for short) over ${\mathbb F}$ is
called a pseudo-Riemannian Lie algebra if it satisfies (\ref{P1}),
(\ref{P2}) and
\begin{equation}\label{P3}
[A_uv,w]+[u,A_wv]=0, \forall u,v,w\in\fg.
\end{equation}

$\fg$ is called a Riemann-Lie algebra if $(\ ,\ )$ is positive
definite.
\end{Def}

For convenience, we define a product on $\g$ by
$$uv=A_uv, \forall u,v\in\g.$$
Then $\g$ has a new algebraic structure. And we can rewrite the
above definition.

\begin{Def}
A pair $(\fg,(\ ,\ ))$ is called a pseudo-Riemannian Lie algebra
if for any $u$, $v$, $w\in\g$, it satisfies the following
conditions:

\hspace{1in}(PR1)\ \  $uv-vu=[u,v]$;

\hspace{1in}(PR2)\ \  $(uv,w)+(v,uw)=0$;

\hspace{1in}(PR3)\ \ $[uv,w]+[u,wv]=0$.
\end{Def}

By (PR1) we can write condition (PR3) in another form\vspace{3mm}

\hspace{1in}({\it PR3'})\ \
$(uv)w-w(uv)+u(wv)-(wv)u=0$.\vspace{3mm}

So we can redefine pseudo-Riemannian Lie algebra from another
point of view.

\begin{Def}
An algebra $A$ with a non-degenerate symmetric bilinear form $(\
,\ )$ is called a pseudo-Riemannian Lie algebra if for any $u$,
$v$, $w\in A$, the conditions (PR2) and (PR3') are satisfied.
\end{Def}

Under this definition, if we define $[u,v]=uv-vu$, which defines a
Lie algebra structure on $\g$ since one can easily see that (PR3')
implies the Jacobi identity. Thus $A$ is a Lie-admissible algebra.
So the two definitions are equivalent. Denote by $l_u=A_u$ and
$r_u$, $u\in\g$, the left and right multiplication by $u$,
respectively. Then (PR2) and (PR3) can be written as
$$(l_uv,w)+(v,l_uw)=0,\ \ \ [r_vu,w]+[u,r_vw]=0.$$ %\vspace{3mm}

\begin{rem} If $\g$ is an abelian Lie algebra, then the product is
trivial, i.e., $xy=0$, for any $x$, $y\in\g$.
\end{rem}

Now we can state our main theorem in this paper.

\begin{thm}\label{main}
Assume that $\F$ is a field of characteristic 0. Then any Lie
algebra over $\F$ with a product satisfying (PR1) and (PR3) is
solvable. Consequently any pseudo-Riemannian Lie algebra $(\g,(\
,\ ))$ over $\F$ is solvable.
\end{thm}

\section{The Proof of Theorem~\ref{main}}\label{secm}

Firstly, we prove the following lemma, which is a decisive result
leading to our main result.

\begin{lem}
Let $\fg$ be a Lie algebra over an algebraically closed field of
characteristic 0 with a product satisfying (PR1) and (PR3). Then
$\fg$ is not semi-simple.
\end{lem}

\pf Assume that $\fg$ is semi-simple and $\h$ is a Cartan
subalgebra of $\g$. Let $\Delta(\g,\h)$ be the root system of $\g$
with respect to $\h$. Choosing any system of positive roots, we
have the root subspace decomposition of $\g$:
$$\fg=\fh+\sum_{\alpha>0}\fg_\alpha+\sum_{\alpha<0}\fg_\alpha.$$
Let $\fg^+=\sum_{\alpha>0}\fg_\alpha$ and
$\fg^-=\sum_{\alpha<0}\fg_\alpha$. Since $\fg \supseteq \fg\fg
\supseteq [\fg,\fg]$, then $\fg=\fg\fg=[\fg,\fg].$

Let $X_{\alpha}\in\g_{\alpha}$, $\alpha\in\Delta(\g,\fh)$, be a
Chevalley basis of $\g$ with respect to $\h$. We will prove the
lemma in the following steps:
\begin{enumerate}
   \item $\h\h\subset \h.$\\ For any $h_1$, $h_2\in\fh$,
   assume that $$h_1h_2=h_0+X^++X^-,$$ where $h_0\in\fh$, $X^+\in \fg^+$ and $X^-\in
   \fg^-$. If $X^+\neq 0$, then there exists $Y^-\in \fg^-$ such
   that the projection of $[X^+,Y^-]$ to $\fh$ is nonzero.
   By (PR3),
   $$[h_1h_2,Y^-]=-[h_1,Y^-h_2].$$
   But the projection to $\h$ of the left hand side is nonzero while that of
   the right hand side is zero since $[\h,\fg^-]\subset
   \fg^-$, $[\fg^-,\fg^-]\subset \fg^-$, $[\h,\fg^+]\subset
   \fg^+$ and $[\h,\h]=0$. That is a contradiction. Therefore,
   $X^+=0$. Similarly, $X^-=0$. Namely, $\h\h\subset \h.$

   \item $\fg_\alpha \h\subset \fg_\alpha$. Consequently $\h\g_{\alpha}\subset\g_{\alpha}$.\\
   For any $h_1,h_2\in\fh$, one has $$[X_\alpha
   h_1,h_2]+[X_\alpha,h_2h_1]=0.$$ Since $\h\h\subset \h$, then
   $[X_\alpha,h_2h_1]\in \fg_\alpha$, which implies $X_\alpha
   h_1\in \fg_\alpha +\h$ since $h_2$ is arbitrary. Assume that $X_\alpha
   h_1=cX_\alpha +h_0$, where $c\in\C$ and $h_0\in\fh$. If $h_0\neq 0$, then there exists a root $\beta$
   such that $\beta(h_0)\not=0$. We can assume $\beta\neq \alpha$
   since if $\alpha(h_0)\neq 0$, we can choose $\beta=-\alpha.$
   Then, $$[X_\alpha h_1,X_\beta]=[cX_\alpha+h_0,X_\beta]=
   cN_{\alpha,\beta}X_{\alpha+\beta}+\beta(h_0)X_\beta,$$
   where $N_{\alpha,\beta}$ are the Chevalley coefficients.
   Similarly, $$[X_\alpha,X_\beta h_1]=[X_\alpha,c'X_\beta+h'_0]=
   c'N_{\alpha,\beta}X_{\alpha+\beta}-\alpha(h'_0)X_\alpha.$$
   Then, $[X_\alpha h_1,X_\beta]+[X_\alpha,X_\beta h_1]\not=0$, which
   contradicts to ({\it PR3}).
   Thus, $h_0=0$. Namely, $\fg_\alpha \h\subset
   \fg_\alpha$. By ({\it PR1}), one can get
   $\h\g_{\alpha}\subset\g_{\alpha}$.

   \item $X_\alpha h=f(h)X_\alpha$, and
   $X_{-\alpha}h=-f(h)X_{-\alpha}$ for some $f\in\h^*$.

   From the above discussion, we may
   assume $$X_\alpha
   h=f(h)X_\alpha,\ \ \ X_{-\alpha}h=g(h)X_{-\alpha},$$
   for some $f$, $g\in\h^*$, since $\dim \fg_\alpha=1$.
   By ({\it PR3}) $$[X_\alpha h,X_{-\alpha}]+[X_\alpha,X_{-\alpha}
   h]=0,$$namely,
   $$[f(h)X_\alpha,X_{-\alpha}]+[X_\alpha,g(h)X_{-\alpha}]=0.$$
   Then, $(f(h)+g(h))[X_\alpha,X_{-\alpha}]=0$. Therefore, $f(h)+g(h)=0$.

%   \item $\h\fg_\alpha \subset \fg_\alpha.$\\ Since $\fg_\alpha \h\subset \fg_\alpha
%   $ and $hX_\alpha=X_\alpha h+[h,X_\alpha]$, then
%   $hX_\alpha\subset \fg_\alpha$.

   \item $\fg_\alpha \h=0$.\\
   For any root $\alpha$, there exists $h_1\in \h$ such that
   $\alpha(h_1)\not=0$.
   For any $h_2\in\fh$, we
   have $$[h_1h_2,X_\alpha]+[h_1,X_\alpha h_2]=0$$ and $$[h_1h_2,X_{-\alpha}]+[h_1,X_{-\alpha}
   h_2]=0.$$
   Since $\h\h\subset \h$, $X_\alpha h=f(h)X_\alpha$ and $X_{-\alpha}
   h=-f(h)X_\alpha$,  we get
$$\alpha(h_1h_2)X_\alpha+f(h_2)\alpha(h_1)X_\alpha=0;$$
$$(-\alpha)(h_1h_2)X_{-\alpha}+(-f)(h_2)(-\alpha)(h_1)X_{-\alpha}=0.$$
   Thus, $$\alpha(h_1h_2)+f(h_2)\alpha(h_1)=0$$
   $$(-\alpha)(h_1h_2)+f(h_2)\alpha(h_1)=0.$$ Therefore, $f(h_2)=\alpha(h_1h_2)=0$.
   So $X_\alpha h_2=f(h_2)X_{\alpha}=0$. Since $h_2$ is arbitrary, we get $X_\alpha
   \h=0$, that is, $\fg_\alpha \h=0$.

\item $\h\h=0$.

   Assume that $h_1h_2\not=0$, for some $h_1$, $h_2\in\h$.
   Then there exists a root $\alpha$ such that
   $\alpha(h_1h_2)\not=0$. Thus
   $$[h_1h_2,X_\alpha]=\alpha(h_1h_2)X_\alpha\not=0.$$
   But by ({\it PR3}) and $\g_\alpha\fh=0$, we have
   $$[h_1h_2,X_\alpha]=-[h_1,X_\alpha h_2]=0.$$
   That is a contradiction. Therefore, $h_1h_2=0$, i.e., $\h\h=0$.

   \item $\fg_\alpha\fg_\beta\subset \fg_{\alpha+\beta}.$ Here $\g_0=\h$.\\ For
   any $h\in\h$, one has
   $$[h,X_\alpha X_\beta]+[hX_\beta,X_\alpha]=0$$
   But $[hX_\beta,X_\alpha]\in \fg_{\alpha+\beta}$
   since $hX_\beta\in \fg_\beta$. Thus $[h,X_\alpha X_\beta]\in
   \fg_{\alpha+\beta}$, for any $h\in\h$. Therefore we have $$X_\alpha X_\beta\in
   \fg_{\alpha+\beta}+\h.$$ Assume that $X_\alpha X_\beta=c_{\alpha,\beta}
   X_{\alpha+\beta}+h_1$, where $c_{\alpha,\beta}\in\C$ and
   $h_1\in\fh$. If $\alpha+\beta=0$, we've done since
   $\g_{\alpha+\beta}=\g_0=\h$. So in the following, we assume
   $\alpha+\beta\neq 0$.
   If $h_1\not=0$, then there exists a root $\gamma$
   such that $\gamma(h_1)\not=0$. We can assume $\gamma\neq
   \alpha$, since if $\alpha(h_1)\neq 0$, we can choose
   $\gamma=-\alpha$.
   Then $$[X_\gamma,X_\alpha X_\beta]=[X_\gamma,c_{\alpha,\beta}
   X_{\alpha+\beta}+h_1]=c_{\alpha,\beta}[X_\gamma,X_{\alpha+\beta}]-\gamma(h_1)X_\gamma.$$
   Assume that $X_\gamma X_\beta=c_{\gamma,\beta}X_{\gamma+\beta}+h_2$. Then
   $$[X_\gamma X_\beta,X_\alpha]=[c_{\gamma,\beta}
   X_{\gamma+\beta}+h_2,X_\alpha]=c_{\gamma,\beta}[X_{\gamma+\beta},X_{\alpha}]+\alpha(h_2)X_\alpha.$$
   Since $\alpha+\beta\neq 0$, $\alpha\neq\gamma$ and $\gamma(h_1)\neq 0$,
   one has $[X_\gamma,X_\alpha X_\beta]+[X_\gamma
   X_\beta,X_\alpha]\not=0$. That is a contradiction. Thus $h_1=0$.
   So $X_\alpha X_\beta\in
   \fg_{\alpha+\beta}$.

   \item $\h\fg_\alpha=0$.\\
   Assume that $hX_\alpha\neq 0$, for some $h\in\h$. Then $hX_\alpha=f(h)X_\alpha$ and $f(h)\not=0$.
By ({\it PR3}),
   $$[h,X_{-\alpha}X_\alpha]+[hX_\alpha,X_{-\alpha}]=0.$$ Since
   $X_{-\alpha}X_\alpha\in\g_0=\h$, one has
   $$[h,X_{-\alpha}X_\alpha]+[hX_\alpha,X_{-\alpha}]=[hX_\alpha,X_{-\alpha}]=f(h)[X_\alpha,X_{-\alpha}]\not=0.$$
That is a contradiction. Thus $\h\fg_\alpha=0$.
\end{enumerate}

Now there exists a contradiction since $[\h,\fg_\alpha]=0$ by
$\h\g_\alpha=\g_\alpha\h=0$. Thus, $\fg$ is not a semi-simple Lie
algebra.\hfill$\Box$

\vspace{3mm}

Now we come to the proof of our main result.

{\it The proof of the main theorem.}\ \ We can assume that $\F$ is
algebraically closed, otherwise we can extend $\F$ to its
algebraic closure. Let $\fg=\mathfrak s+\mathfrak r$ be a Levi
decomposition of $\fg$. For any $s_1,s_2,s_3\in \mathfrak s$,
$$[s_1s_2,s_3]+[s_1,s_3s_2]=0.$$

Let $s_is_j=s_{i,j}+r_{i,j}$, where $s_{i,j}\in \mathfrak s$ and
$r_{i,j}\in \mathfrak r$. Then
$$[s_{1,2}+r_{1,2},s_3]+[s_1,s_{3,2}+r_{3,2}]=0,$$ i.e.,
$$([s_{1,2},s_3]+[s_1,s_{3,2}])+([r_{1,2},s_3]+[s_1,r_{3,2}])=0.$$
Thus
$$[s_{1,2},s_3]+[s_1,s_{3,2}]=[r_{1,2},s_3]+[s_1,r_{3,2}]=0,$$
since $\mathfrak s$ is a subalgebra and $\mathfrak r$ is an ideal
of $\fg$. Define a product $\circ$ of $\mathfrak s\times \mathfrak
s\rightarrow \mathfrak s$ by $$s_1\circ s_2=P_s(s_1s_2),$$ where
$P_s$ denotes the projection from $\fg$ to $\mathfrak s$ with
respect to the Levi decomposition. Then we have
\begin{enumerate}
   \item The product $\circ$ is bilinear.
   \item For any $s_1,s_2\in \mathfrak s$,
    $$[s_1,s_2]=s_1s_2-s_2s_1=s_{1,2}+r_{1,2}-s_{2,1}-r_{2,1}=(s_{1,2}-s_{2,1})+(r_{1,2}-r_{2,1})\in \mathfrak
   s.$$Thus $r_{1,2}-r_{2,1}=0$. Therefore,
   $$s_1\circ s_2-s_2\circ s_1=P_s(s_1s_2)-P_s(s_2s_1)=s_{1,2}-s_{2,1}=[s_1,s_2].$$
   \item For any $s_1,s_2,s_3\in \mathfrak s$,
        $$[s_1\circ s_2,s_3]+[s_1,s_3\circ s_2]
         =
         [P_s(s_1s_2),s_3]+[s_1,P_s(s_3s_2)]=[s_{1,2},s_3]+[s_1,s_{3,2}]=0.$$
\end{enumerate}
Thus, $(\mathfrak s,\circ)$ satisfy the conditions of the above
lemma. Therefore, $\mathfrak s$ is not semi-simple. Then,
$\mathfrak s$ must be $0$. Namely, $\fg=\mathfrak r$ is
solvable.\hfill$\Box$

\section{Riemann-Lie Algebras}\label{sec4}
In this section, an ideal of a pseudo-Riemannian Lie algebra $\g$
means it is invariant under left and right multiplications in
$\g$, hence an ideal is automatically a Lie-ideal.
% i.e., an ideal as a Lie subalgebra.

\subsection{Ideals of $\g$ as a Lie-admissible Algebra}

Let $C(\g)=\{a\in \g|ax=xa,\forall x\in \g\}$ be the center of
$\g$, then we have
\begin{lem}
$C(\g)$ is an ideal of $\g$, and $xy=0$ for any $x$, $y\in C(\g)$.
\end{lem}
\pf For any $a\in C(\g)$, $x,y\in \g$, one has by $(PR3)$
$$[ax,y]+[a,yx]=[ax,y]=0,$$
which implies that $ax=xa\in C(\g)$, $\forall x\in \g$.

For any $x$, $y\in C(\g)$ and $z\in\g$, then
$$(xy,z)=-(y,xz)=-(y,zx)=(zy,x)=(yz,x)=-(z,yx)=-(xy,z),$$
which implies $(xy,z)=0$, for any $z\in\g$. Thus $xy=0$.
\hfill$\Box$

\begin{rem} This lemma is interesting and non-trivial
since the center of general algebras, say, associative algebras,
is not necessary an ideal.
\end{rem}
%\vspace{3mm}

For any $x\in C(\g)$, $y\in C(\g)^{\bot}=\{u\in\g|(u,C(\g)=0\}$,
$z\in\g$, we have
$$(x,yz)=-(yx,z)=-(xy,z)=(y,xz)=0=(y,zx)=(zy,x).$$
So $yz$ and $zy\in C(\g)^{\bot}$. Thus we have

\begin{lem}\label{lem4.2}
$C(\g)^{\bot}$ is an ideal of $\g$. If the restriction of the
bilinear form to $C(\g)$ is non-degenerate (say, if $\g$ is a
Riemann-Lie algebra), then $[\g,\g]\subset C(\g)^{\bot}$ and
$$\g=C(\g)\oplus C(\g)^{\bot}.$$
\end{lem}

%Let $Z(\g)=\{z\in \g|zx=xz=0,\forall x\in \g\}$, which is an
%abelian ideal of $\g$ contained in $C(\g)$, and  $Z_l(\g)=\{z\in
%\g|zx=0,\forall x\in \g\}$ and $Z_r(\g)=\{z\in \g|xz=0,\forall
%x\in \g\}$.

set $\g\g=<\{xy|x,y\in\g\}>$ and $Z_r(\g)=\{u\in \g|r_u=0\}$. Then
we have
\begin{lem}\label{lem4.3}
$\g\g$ is an ideal of $\g$ and $(\g\g)^\bot=Z_r(\g)$.
\end{lem}
\pf The first assertion is trivial. For the second one, $\forall
x,y,z\in\g$, $$x\in(\g\g)^\bot\Leftrightarrow(x,yz)=0,
\Leftrightarrow(yx,z)=0\Leftrightarrow yx=0\Leftrightarrow x\in
Z_r(\g).\eqno{\Box}$$

\subsection{New Proof of M. Boucetta's Results}

In this subsection, we will use our new definition and result to
give a short and simple proof of M. Boucetta's results in
\cite{MB3}. An immediate consequence of our main theorem is that
Lemma~3.5 in~\cite{MB3} is trivial. %We will try to approach the
%results under the assumption that $\g$ is a pseudo-Riemannian Lie
%algebra and avoid using that $\g$ is Riemann-Lie algebra until
%necessary.

Let $\g$ be a Lie algebra with a scalar product $(\ ,\ )$. Define
the adjoint $\phi^t$ of $\phi\in End(\g)$ by
$$(\phi(v),w)=(v,\phi^t(w)),\ \ \ \forall v,w\in\g.$$

First, we have the following Lemmas.

\begin{lem} $[\g,\g]^\bot=\{u\in\g|r_u=r_u^t\}$. %For $u$,
%$v\in[\g,\g]^\bot$, $uv+vu=0$.
\end{lem}
\pf For any $u$, $v$, $w\in\g$,
$$u\in[\g,\g]^{\bot}\Leftrightarrow(u,[v,w])=0\Leftrightarrow(u,vw)=(u,wv)\Leftrightarrow
(vu,w)=(v,wu)\Leftrightarrow r_u=r_u^t.$$ So
$[\g,\g]^\bot=\{u\in\g|r_u=r_u^t\}$.\hfill$\Box$\vspace{3mm}

%For $u$, $v\in [\g,\g]^\bot$, $w\in\g$, we have
%$$(uv,w)=(u,wv)=-(wu,v)=-(w,vu).$$
%So $(uv+vu,w)=0$, i.e., $uv+vu=0$.\hfill$\Box$%\vspace{3mm}

Notice that for any $u\in[\g,\g]^\bot$, $r_u$ is self-adjoint, so
if the bilinear form is positive definite, then $r_u$ is
diagonalizable. Then one can easily prove
\begin{lem}[\cite{MB3}, Lemma~3.3] Let $(\g,(\ ,\ ))$ be a Riemann-Lie algebra. Then
$[\g,\g]^\bot=Z_r(\g)$. \hfill$\Box$
\end{lem}

By Lemma~\ref{lem4.3}, $[\g,\g]=\g\g$ provided $\g$ is a
Riemann-Lie algebra, thus $[\g,\g]$ is a Riemann-Lie subalgebra.
Hence $[\g,\g]$ is abelian by Lemma~\ref{lem4.2} since $[\g,\g]$
is nilpotent. Thus we have the decomposition of $\g$ into an
orthogonal direct sum
$$\g=Z_r(\g)\dot+[\g,\g]$$
with $Z_r(\g)$ an abelian subalgebra and $[\g,\g]$ an abelian
ideal. So we have almost proved the following result of
M.Boucetta.

\begin{thm}[\cite{MB3}, Theorem 2.1]
Let $\g$ be a real Lie algebra and $(\ ,\ )$ a scalar product on
$\g$. Then the following assertions are equivalent:

1) $(\g,(\ ,\ ))$ is a Riemann-Lie algebra.

2) The orthogonal Lie-subalgebra $S$ of $(\g,(\ ,\ ))$ is abelian
and $\g$ splits as an orthogonal direct sum $S\dot + U$, where $U$
is a commutative Lie-ideal and $S=\{s\in\g|\ad_s+\ad_s^t=0\}$,
which is a Lie-subalgebra.
\end{thm}

In order to complete the proof of the above theorem, we just need
to show $S=Z_r(\g)$, or equivalently
$S^\bot=\{u\in\g|(u,S)=0\}=\g\g$.

\begin{lem}
$S^\bot=\{uv+vu|u,v\in\g\}=\g\g.$
\end{lem}
\pf %For any $s$, $x\in S$, $u\in S^\bot$, we have $sx=0$, so
%$0=(sx,u)=-(x,su)=-(s,ux)$. Hence $su$, $ux\in S^\bot$.
For any $u$, $v\in\g$, $s\in S$, we have
$$(s,uv)=-(us,v)=(u,vs)=-(uv,s),$$
which implies $(s,uv+vu)=0$. Notice that $\g=\g\g\dot+ Z_r(\g)$
and $\g\g$ and $Z_r(\g)$ are abelian, then
$\g\g=\{uv+vu|u,v\in\g\}$. So we have
$S^\bot=\g\g$.\hfill$\Box$\vspace{3mm}

\subsection{Structure and Construction of Riemann-Lie Algebras}

From the above discussion, one can easily know the structure of
$\g$: the product and the bilinear form, and can easily construct
some Riemann-Lie algebras.

\begin{cor}\label{cor4.11}
Let $\g$ be a 3-dimensional non-abelian Riemann-Lie algebra. Then
$\dim[\g,\g]=2$. There exists an orthonormal basis $s$, $x$, $y$
of $\g$ such that $s\in Z_r(\g)$ and $x$, $y\in[\g,\g]$, and
$[s,x]=sx=ay$, $[s,y]=sy=-ax$, for some $a\in\R$.
\end{cor}

\begin{rem}
If we define $<u,v>=a^2(u,v)$, then $s'=\frac{1}{a}s$,
$x'=\frac{1}{a}x$, $y'=\frac{1}{a}y$ form an orthonormal basis of
$\g$, and $[s',x']=s'x'=y'$, $[s',y']=s'y'=-x'$. In other words,
there is a unique inner product on the above Lie algebra $\g$ (up
to a positive constant) such that $\g$ is a Riemann-Lie algebra.
\end{rem}

Now $\g=Z_r(\g)\dot+\g\g$, with $\g\g=[\g,\g]$. Furthermore, we
have

1) $us=0$, $su=[s,u]$ for all $u\in\g$, $s\in Z_r(\g)$.

2) $xy=0$ for all $x$, $y\in\g\g$.

We can assume $\g$ has no center. Now consider the adjoint
representation $\rho:Z_r(\g)\rightarrow \fso(\g\g)$ by
$$\rho(s)=\ad_s|_{\g\g}.$$
One can see the scalar product on $\g\g$ is a $Z_r(\g)$-invariant
inner product. So the representation is completely reducible and
it is decomposed into irreducible sub-representations of
2-dimensional (no 1-dimensional (trivial) sub-representation!)
$$\g\g=L_1\dot+\cdots\dot+ L_n.$$
Assume that $Z_r(\g)=H_1\dot+\cdots\dot+ H_m$ is an orthogonal
decomposition of $Z_r(\g)$ into 1-dimensional subspaces. Then one
can see $H_i\dot+ L_j$, $1\leq i\leq m$, $1\leq j\leq n$ are
Riemann-Lie subalgebras. So one can construct Riemann-Lie algebras
in the following way. First let $H\dot+ L_1,\ldots,H\dot+ L_n$ be
3-dimensional Riemann-Lie algebras, then $H\dot+ L=H\dot+
L_1\dot+\cdots\dot+ L_n$ is a Riemann-Lie algebra (with the
naturally algebraic structure). Choosing Riemann-Lie algebras
$H_1\dot+ L,\ldots,H_m\dot+ L$, one can see
$\g=(H_1\dot+\cdots\dot+ H_m)\dot+ L$ is also a Riemann-Lie
algebra.

\section{Classification of Pseudo-Riemannian Lie Algebras of Dimension $\leq
3$}\label{sec5} In this section we classify the pseudo-Riemannian
Lie algebras of dimension $\leq 3$. Let $(\g,(\ ,\ ))$ be a
non-abelian pseudo-Riemannian Lie algebra. Since $\g$ is solvable,
$\dim C(\g)\leq 1$. M. Boucetta claimed the same classification in
Theorem~1.6~\cite{MB2} without proof. By our definition and
method, the proof should be easier than his. Furthermore,
Theorem~1.6 2b) is not quite correct and M. Boucetta didn't give
the concrete formula and bilinear form for $\g$ being a
pseudo-Riemannian Lie algebra, which we describe explicitly in the
following. Anyway, the method is different.

\subsection{$\dim\g=2$}

Since $\g$ is non-abelian, then $\g$ has a basis $x,y$ satisfying
$[x,y]=y$.

If $(y,y)=0$, then $(x,y)\neq 0$. Replacing $x$ by
$x-\frac{(x,x)}{2(x,y)}y$, one can assume that $(x,x)=0$. Since
$(xy,y)=0$, then $xy\in<y>$. $(xx,x)=(yx,x)=0$ implies
$xx,yx\in<x>$. Then $[x,y]=y$ implies that $xy=y$, $yx=0$,
%Similarly, we can assume $xx=ax$, $yy=by$ for some $a,b\in\F$.
and $[xx,y]+[x,yx]=0$ implies $xx=0$. So $(xy,x)+(y,xx)=0$ implies
$(y,x)=0$. That is a contradiction.

If $(y,y)\neq 0$, similar argument deduces a contradiction. So we
have
%replacing $x$ by $x-\frac{(x,y)}{(y,y)}y$, one can
%assume $(x,y)=0$, and $(x,x)\neq 0$.

%Since $(xy,y)=0$, then $xy=kx$, for some $k\in\F$. $(yx,x)=0$
%implies $yx=ly$ for some $l\in\F$. But $[x,y]=y$ implies that
%$k=0$, $l=-1$. Similarly, we can assume $xx=ay$, $yy=bx$ for some
%$a,b\in\F$. But $[xx,y]+[x,yx]=0$ implies $[x,yx]=[x,-y]=0$.
%Contradiction.
\begin{lem}
2-dimensional Riemann-Lie algebra  must be abelian.
\end{lem}

\subsection{$\dim\g=3$: $\dim C(\g)=1$}

In this case, $\g$ is a Heisenberg Lie algebra or a direct sum of
$C(\g)$ and the two dimensional non-abelian Lie algebra.

\begin{prop}
The only 3-dimensional non-abelian pseudo-Riemannian Lie algebra
$\g$ with non-trivial center is a Heisenberg Lie algebra. And
assume that $x$, $y$, $z$ is a standard basis of $\g$ such that
$[x,y]=z$, then the product and the scalar product can be given as
follows.

1) $(x,z)=1$, $(y,y)\neq 0$ and otherwise zero.

2) $xx=-\frac{1}{(y,y)}y$, $xy=z$ and otherwise zero.

Furthermore, $\g$ can't be a Riemann-Lie algebra.
\end{prop}
\pf Assume that $z$ is a basis of $C(\g)$. Then we have $z\g=\g
z\subset C(\g)$, $zz=0$ and $(z,z)=0$. Thus $C(\g)\subset
C(\g)^\bot$. Assume that $y$, $z$ is a basis of $C(\g)^\bot$, then
$(y,y)\neq 0$. Choose $x\in\g$ such that $(x,x)=(x,y)=0$ and
$(x,z)\neq 0$. We can assume $(x,z)=1$. Then $(zx,x)=0$ implies
$zx=xz=0$ and $(zy,x)=-(y,zx)=0$ implies $zy=yz=0$. So
$\g\g\subset C(\g)^\bot=<y,z>$. Thus we have (notice all the
products are in $<y,z>$)

1) $xx\in<y>$ since $(xx,x)=0$.

2) $yx=0$ since $[x,yx]=-[xx,y]=0$ and $(yx,x)=0$.

3) $xy\in<z>$ since $(xy,y)=0$. So $[x,y]=xy\in<z>$, which means
$\g$ is a Heisenberg Lie algebra. We can assume $[x,y]=z$.

4) $yy=0$ since $(yy,x)=-(y,yx)=0$ and $(yy,y)=0$.

5) $(xx,y)=-(x,xy)=-1$, so $a=-\frac{1}{(y,y)}$. \qed
\begin{rem}
In this case, we have $\g\g\neq[\g,\g]$.
\end{rem}

\subsection{Centerless 3-dimensional Lie Algebras}

Let $(\g,(\ ,\ ))$ be a 3-dimensional centerless pseudo-Riemannian
Lie algebra. Since $\g$ is solvable, there exists a basis $x$,
$y$, $z$ of $\g$ such that
$$[x,y]=ay+bz,\ \ [x,z]=cy+dz,\ \ [y,z]=0,$$
where $ad-bc\neq 0$.
\begin{prop}
Let $\g$ be a Lie algebra over a field $\F$ with a non-degenerate
symmetric bilinear form $(\ ,\ )$. Then the following assertions
are equivalent.

a) $\g$ is pseudo-Riemannian Lie algebra.

b) There exists a basis $x$, $y$, $z$ of $\g$ such $[x,y]=z$,
$[x,z]=cy$ with $c\neq 0$.

In this case, the product can be given by $xy=[x,y]$, $xz=[x,z]$,
otherwise 0. And the bilinear form is given by $(x,x)=a$,
$(y,y)=b$, $(z,z)=-bc$ for any $a\neq 0$, $b\neq 0$, otherwise 0.

If $\F=\R$ and $c<0$, choose $a>0$, $b>0$, one can see $\g$ is a
Riemann-Lie algebra. Choose an orthonormal basis, we can get
Corollary~\ref{cor4.11}.
\end{prop}
\pf First we assume the bilinear form restricting to $<y,z>$ is
degenerate. Thus we can assume $z\in <y,z>^\bot$. So $(y,y)\neq
0$. Then we can choose $x$ such that $(x,y)=0$, $(x,x)=0$ and
$(x,z)=1$. Then we have the following conclusions:

1) $yz$, $zy$, $yy$, $zz\in<z>$.

By (PR2), $(yz,z)=0=(zy,y)$, so $zy=yz\in <y,z>^\bot=<z>$.
Similarly one can prove $yy$, $zz\in<z>$.

2) $xz$, $xy\in<y,z>$ and $yz=0$.

Since $(xz,z)=0\Rightarrow xz\in<y,z>$ and $zx\in<y,z>$, then
$0=[xz,y]=-[x,yz]=-[x,zy]=[xy,z]$, which implies
$$yz=0,\ \ \ xy\in<y,z>.$$
%Furthermore, by $[x,yy]=-[xy,y]=0$ we have $yy=0$.

3) $zx=0$.

$(zx,y)=-(x,zy)=0$ and $<zx,x>=0$ imply that $zx=0$ since
$zx\in<y,z>$.

4) $xx\in<y,z>$.

$[xx,z]=-[x,zx]=0$, so $xx\in<y,z>$.

The above results show that $\g\g\subset[\g,\g]$, hence
$\g\g=[\g,\g]$. So $Z_r(\g)=(\g\g)^\bot=[\g,\g]^\bot=<z>$. So
$[x,z]=xz=0$, which contradicts to the fact that $\g$ has no
center.

So the restriction to $[\g,\g]$ of the bilinear form is
non-degenerate, then we can choose $x$, $y$, $z$ orthogonal.
Notice that $r_x^t=r_x$ since $x\in[\g,\g]^\bot$.

For any $u\in\g$, since $(u,xx)=(ux,x)=0$, then $ux\in<y,z>$ and
$xx=0$, so $[ux,x]=-[u,xx]=0$, which implies
$ux\in\Ker\ad_x\cap<y,z>=\{0\}$, i.e., $x\in Z_r(\g)$. Then
$\g\g=[\g,\g]$. Thus by $(yz,z)=0=(zy,y)$, we have $yz=zy=0$.
Similarly, $yy=zz=0$. Furthermore $(xy,y)=0$ implies
$[x,y]=xy\in<z>$, so $a=0$. Similarly, $d=0$. By $(xz,y)=-(z,xy)$,
one has $c(y,y)=-b(z,z)$. Replacing $x$ by $\frac{1}{b}x$, we
complete the proof.\qed

\end{document}